\documentclass[12pt]{scrartcl}

\usepackage[english]{babel}
\RequirePackage[T1]{fontenc}
\RequirePackage[utf8]{inputenc}

\RequirePackage{geometry}
\RequirePackage{palatino}
\RequirePackage[osf,sc]{mathpazo} 
\RequirePackage{mathrsfs} 
\RequirePackage{stmaryrd} 
\RequirePackage{microtype} 

\RequirePackage[leqno]{amsmath}
\RequirePackage{amsfonts}
\RequirePackage{amssymb}
\usepackage[svgnames]{xcolor} 
\RequirePackage{graphicx}

\RequirePackage{pgf}
\RequirePackage{tikz}

\RequirePackage{marginnote}  
\usepackage[alwaysadjust]{paralist}
\setdefaultenum{\footnotesize (1.)}{\footnotesize (a)}{\footnotesize (i.)}{\footnotesize A.}

\RequirePackage[pdftex,pagebackref,plainpages=false, pdfpagelabels,unicode]{hyperref}

\newenvironment{prose}{\begin{quotation}\footnotesize\noindent[[}
{]]\end{quotation}}

\providecommand{\Marginnote}[2][0pt]{\marginnote{\small {\color{RoyalBlue}#2}}[#1]}

\RequirePackage[amsmath,thmmarks,thref,hyperref]{ntheorem}

{
\theoremstyle{break}
\newtheorem{preTheorem}{Theorem}[section]
\newtheorem{preLemma}[preTheorem]{Lemma}
\newtheorem{preCorollary}[preTheorem]{Corollary}
\newtheorem{preProposition}[preTheorem]{Proposition}
}
{\theoremstyle{plain}
\theorembodyfont{\normalfont}
\newtheorem{preDefinition}[preTheorem]{Definition}
\newtheorem{preRemark}[preTheorem]{Remark}
\newtheorem{preExample}[preTheorem]{Example}
\newtheorem{preQuestion}[preTheorem]{Question}
}

\newenvironment{Theorem}{\bgroup \begin{preTheorem} }{\end{preTheorem} \egroup}
\newenvironment{Lemma}{\bgroup \begin{preLemma} }{\end{preLemma} \egroup}
\newenvironment{Corollary}{\bgroup \begin{preCorollary} }{\end{preCorollary} \egroup}
\newenvironment{Proposition}{\bgroup \begin{preProposition}}{\end{preProposition} \egroup}

\newenvironment{Definition}{\bgroup \begin{preDefinition} }{\end{preDefinition} \egroup}
\newenvironment{Remark}{\bgroup \begin{preRemark} }{\end{preRemark} \egroup}

\newenvironment{Question}{\bgroup \begin{preQuestion} }{\end{preQuestion} \egroup}

\newcommand{\plist}{\begin{compactenum}}
\newcommand{\pliste}{\end{compactenum}}

{
\theoremstyle{nonumberplain}
\theoremheaderfont{\normalfont\slshape}
\theorembodyfont{\upshape}
\theoremseparator{.}
\theoremsymbol{\ensuremath{_\Box}}
\newtheorem{proof}{Proof}
}

{
\theoremstyle{nonumberplain}
\theoremheaderfont{\normalfont}
\theorembodyfont{\upshape \footnotesize}
\theoremseparator{}
\theoremsymbol{ \ensuremath{_\boxtimes}}

\newtheorem{spoiler}{Summary.}
}

\newcommand{\hl}[1]{{\color{RoyalBlue}#1}}

\mathchardef\mathhyphen="2D 

\newcommand{\widebar}[1]{\overline{#1}}
\newcommand{\bbar}{\widebar} 

\newcommand{\matn}{\mathbb{N}}
\newcommand{\matz}{\mathbb{Z}}

\newcommand{\matf}{\mathbb{F}}

\newcommand{\mathp}{\mathfrak{P}}

\newcommand{\restr}{\!\upharpoonright\!}
\newcommand{\notnil}{\not= \emptyset}
\newcommand{\nil}{\emptyset}
\newcommand{\phii}{\varphi}

\newcommand{\mal}{\cdot}

\newcommand{\xsupp}{\mathrm{\mathhyphen supp}}
\newcommand{\xmin}{\mathrm{\mathhyphen min}}
\newcommand{\xmax}{\mathrm{\mathhyphen max}}
\newcommand{\ZFC}{\textsl{ZFC}}
\newcommand{\CH}{\textsl{CH}}

\newcommand{\dotcup}{\dot{\cup}}

\hypersetup{
   pdftitle={On union ultrafilters},    
   pdfauthor={Peter Krautzberger},     
   pdfsubject={},   
   pdfkeywords={ mathematics, ultrafilters, stone-cech compactification, combinatorics, set theory, idempotent filter,forcing, idempotent ultrafilter,union ultrafilter, summable ultrafilter, strongly summable ultrafilter}, 
   linkbordercolor = RoyalBlue,
   citebordercolor = RoyalBlue,
   urlbordercolor = RoyalBlue}

\title{On union ultrafilters}
\author{Peter Krautzberger\thanks{Mathematics Department, University of Michigan, Ann Arbor, pkrautzb@umich.edu, Partially supported by DFG-grant KR 3818; Subject classification 03E75 (Primary) 54D80, 05D10 (Secondary)} }
\date{\today}

\begin{document}

\maketitle
\begin{abstract}
We present some new results on union ultrafilters. We characterize stability for union ultrafilters and, as the main result, we construct a new kind of unordered union ultrafilter.
\end{abstract}

\section*{Introduction}

The equivalent notions of union and strongly summable ultrafilters have been important examples of idempotent ultrafilters ever since they were first conceived in \cite{Hindman72}, \cite{Blass87-1}. Their unique properties have been a useful tool in set theory, algebra in the Stone-\v Cech compactification and set theoretic topology. For example, strongly summable ultrafilters were, in a manner of speaking, the first idempotent ultrafilters known, cf.~\cite{Hindman72} and \cite[notes to Chapter 5]{HindmanStrauss}; they were the first strongly right maximal idempotents known and, even stronger, they are the only known idempotent with a maximal group isomorphic to $\matz$; their existence is independent of \ZFC, since it implies the existence of (rapid) $P$-points, cf.~\cite{BlassHindman87}; since a strongly summable is strongly right maximal, its orbit is a van Douwen space, cf.~\cite{HindmanStrauss02}.

This article will focus on union ultrafilters, studying the various kinds of union ultrafilters and as the main result constructing of a new kind of union ultrafilter answering a question of Andreas Blass.

The presentation of the proofs is inspired by \cite{Leron83} and \cite{Lamport95} splitting the proofs into different levels, at times adding {\footnotesize [[in the elevator]]} comments in between. The typesetting incorporates ideas from \cite{Tufte05} highlighting details in the proofs and structural remarks in the margin. Online discussion is possible through the author's website at \url{http://peter.krautzberger.info/papers}.

\section{Preliminaries}

Let us begin by giving a non-exhaustive selection of standard terminology in which we follow N.~Hindman and D.~Strauss \cite{HindmanStrauss}; for standard set theoretic notation we refer
to T.~Jech \cite{Jech}, e.g., natural numbers are considered as ordinals, i.e., $n = \{ 0, \ldots , n-1\}$. We work in \ZFC\ throughout. The main objects of this paper are \emph{(ultra)filters} on an infinite set $S$, i.e., (maximal) proper subsets of the power set $\mathp(S)$ closed under taking finite intersections and supersets.  
$S$ carries the discrete topology in which case the set of ultrafilters is $\beta S$, its Stone-\v Cech compactification. The Stone topology on $\beta S$ is generated by basic clopen sets induced by subsets $A\subseteq S$ in the form $ \bbar{A} := \{ p\in \beta S\ |\ A\in p\}$.
Filters are usually denoted by upper case Roman letters, mostly $F,G,H$, ultrafilters by lower case Roman letters, mostly $p,q,r,u$.

The set $S$ is always assumed to be the domain of a \Marginnote{(Partial) Semigroup}
\emph{(partial) semigroup} $(S,\cdot)$, i.e., the (partial) operation $\mal$ fulfills the associativity law  $ s\cdot (t \cdot v) = (s\cdot t) \cdot v $ (in the sense that if one side is defined, then so is the other and they are equal). For a partial semigroup $S$ and $s\in S$ the set of elements compatible with $s$ is denoted by $\sigma(s) := \{ t\in S\ |\ s\cdot t \mbox{ is defined} \}$.
A partial semigroup is also assumed to be \emph{adequate}, i.e., $\{ \sigma(s) \ |\ s\in S\}$ has the finite intersection property.
We denote the generated filter by $\sigma(S)$ and the corresponding closed subset of $\beta S$ by \index{$\delta S$}\index{semigroup!$\delta S$} $\delta S$.
For partial semigroups $S,T$ a map $\phii: S \rightarrow T$ is a \emph{partial semigroup homomorphism}\index{homomorphism!partial semigroup} if $\phii[\sigma(s)] \subseteq \sigma(\phii(s))$ (for $s\in S$) and
\[
(\forall s\in S)(\forall s' \in \sigma(s))\ \phii(s\mal s') = \phii(s) \mal \phii(s').
\]
To simplify notation in a partial semigroup, $s \cdot t$ is always meant to imply $t \in \sigma(s)$. For $s\in S$, the restricted multiplication to $s$ from the left (right) is denoted by $\lambda_s$ ($\rho_s$).

It is easy to see that the operation of a partial semigroup can always be extended to a full semigroup operation by adjoining a (multiplicative) zero which takes the value of all undefined products. One key advantage of partial semigroups is that partial subsemigroups are usually much more diverse than subsemigroups.  Nevertheless, it is convenient to think about most theoretical aspects (such as extension to $\beta S$) with a full operation in mind.

The semigroups considered in this paper are $( \matn , + ) $ (with $\matn := \omega \setminus \{ 0 \}$) and the most important adequate partial semigroup $\matf$.

\begin{Definition}\label{def:matf}\Marginnote{ \normalsize The partial semigroup $\mathbb{F}$}
On \index{semigroup!$\matf$} \index{$\matf$} $\matf := \{ s\subseteq \omega\ |\ \nil \neq s \mbox{ finite} \}$
we define a partial semigroup structure by
\[
s \mal t := s \cup t \mbox{ if and only if } s \cap t = \nil.
\]
\end{Definition}
The theory of the Stone-\v Cech compactification allows for the (somewhat unique) extension of any operation on $S$ to its compactification, in particular a semigroup operation.
\begin{Definition}\label{def:multiplication in beta S}\index{$\beta S$!the semigroup}\Marginnote{\normalsize The semigroup $\beta S$}
For a semigroup $(S,\cdot)$, $s\in S$ and $A \subseteq S$, $p,q \in \beta S$ we define the following.\begin{itemize}
\item $ s^{-1}A := \{ t\in S\ |\ st \in A \} $.
\item $ A^{-q} := \{ s\in S\ |\ s^{-1}A\in q\} $.
\item  $ p \cdot q : = \{ A \subseteq S \ |\ A^{-q} \in p \} $. \\
Equivalently, $p\cdot q$ is generated by sets $\bigcup_{v\in V} v\cdot W_v$ for $V\in p$ and each $W_v  \in q$.
\item $ A^\star := A^{-q} \cap A$. \\
This notation will only be used when there is no confusion regarding the chosen ultrafilter.
\end{itemize}
\end{Definition}
As is well known, this multiplication on $\beta S$ is well defined and extends the operation on $S$. It is associative and right topological, i.e., the operation with fixed right hand side, $\rho_q$, is continuous. For these and all other theoretical background we refer to \cite{HindmanStrauss}.

In the case of a partial semigroup, ultrafilters in $\delta S$ in a way multiply as if the partial operation was total. With the arguments from the following proposition it is a simple but useful exercise to check that if $(S, \cdot)$ is partial the above definitions still work just as well in the sense that $s^{-1}A := \{ t \in \sigma(s) \ | \ st \in A\}$ and $p\cdot q$ is only defined if it is an ultrafilter.
 
\begin{Proposition}\label{prop:semigroup delta S}\Marginnote{\normalsize The semigroup $\delta S$}
Let $S$ be a partial subsemigroup of a semigroup $T$. Then $\delta S$ is a subsemigroup of $\beta T$.
\end{Proposition}
\begin{proof} \plist
\item Simply observe that, by strong associativity, for all $a \in S$
\[
\bigcup_{b \in \sigma(a)} b \cdot ( \sigma ( ab ) \cap \sigma(b) ) \subseteq \sigma(a).
\]
\item Therefore $ \sigma(S) \subseteq p \cdot q$ whenever $p,q \in \delta S$.
\pliste \end{proof} 
It is easy to similarly check that partial semigroup homomorphisms extend to full semigroup homomorphisms on $\delta S$.

Since $A^{-q}$ is not an established notation, the following useful observations present a good opportunity to test it.

\begin{Proposition}\Marginnote{\normalsize Tricks with $A^{-q}$}
Let $p,q \in \beta S$, $A\subseteq S$ and $s, t\in S$.
\begin{itemize}
\item $t^{-1}s^{-1} A= (st)^{-1} A$.
\item $s^{-1}A^{-q} = (s^{-1}A)^{-q}$.
\item $(A \cap B)^{-q} = A^{-q} \cap B^{-q}$.
\item $( s^{-1} A )^\star = s^{-1} A^\star$ (with respect to the same ultrafilter).
\item $( A^{ - q } )^{ - p } = A^{ - ( p \cdot q ) }$.
\end{itemize}
\end{Proposition}
\begin{proof} \plist
This is straightforward to check.
\pliste \end{proof} 

The proverbial big bang for the theory of ultrafilters on semigroups is the following theorem.

\begin{Theorem}[Ellis-Numakura Lemma] \label{thm:ellis-numakura} \index{Ellis-Numakura Lemma}\index{Lemma!Ellis-Numakura}
If $(S,\cdot)$ is a compact, right topological semigroup then there exists an \emph{idempotent} element in $S$, i.e., an element $p \in S$ such that $p\cdot p = p$.
\end{Theorem}
\begin{proof}
See, e.g., \cite[notes to Chapter 2]{HindmanStrauss}.
\end{proof} 

Therefore the following classical fact is meaningful.

\begin{Lemma}[Galvin Fixpoint Lemma]\label{lem:galvin}\index{lemma!{Galvin Fixpoint}}
For idempotent $p\in \beta S$, $A\in p$ implies $A^\star \in p$ and $(A^\star)^\star = A^\star$.
\end{Lemma}

\begin{proof}
$(A^\star)^\star = A^\star \cap (A^\star)^{-p} = A^\star \cap (A \cap A^{-p})^{-p} = A^\star \cap A^{-p} \cap A^{- p\cdot p} = A^\star \cap A^{-p} = A^\star$. 
\end{proof} 

The following definitions are central in what follows. Even though we mostly work in $\mathbb{N}$ and $\mathbb{F}$ we formulate them for a general setting.

\begin{Definition}\label{def:FP-sets etc} \Marginnote{\normalsize $FP$-sets, $\mathbf{x}$-support and condensations} \index{$FP$-set}
Let $\mathbf{x} = (x_n)_{n < N}$ (with $N \leq \omega $) be a sequence in a partial semigroup $(S,\mal)$ and let $K\leq \omega$.
\begin{itemize} 
\item The set of finite products (the \emph{$FP$-set}) is defined as 
\[
FP( \mathbf{x} ) := \{ \prod_{i\in v} x_i\ |\ v\in \matf \},
\]
where products are in increasing order of the indices. In this case, all products are assumed to be defined.\footnote{Note that we will mostly deal with commutative semigroups so the order of indices is not too important in what follows.}

\item $\mathbf{x}$ has \emph{unique representations} if for $v,w \in \matf$ the fact $\prod_{i \in v} x_i = \prod_{j \in w} x_j$ implies $v = w$.

\item If $\mathbf{x}$ has unique representations and $z \in FP(\mathbf{x})$ we can define the \emph{$\mathbf{x}$-support of $z$}, short $\mathbf{x}\xsupp(z)$, \index{$\mathbf{x}\xsupp$} by the equation $z = \prod_{j \in \mathbf{x}\xsupp(z)} x_j$. We can then also define $\mathbf{x}\xmin := \min \circ \mathbf{x}\xsupp$, $\mathbf{x}\xmax := \max \circ \mathbf{x}\xsupp$.


\item A sequence $\mathbf{y}=(y_j)_{j<K}$ is called a \emph{condensation}\index{condensation} of $\mathbf{x }$, in short $\mathbf{y} \sqsubseteq  \mathbf{ x } $, if
\[
FP(\mathbf{y}) \subseteq  FP( \mathbf{x} ).
\]
In particular, $\{ y_i \ |\ i<K \} \subseteq FP(\mathbf{x})$. For convenience, $\mathbf{x}\xsupp(\mathbf{y}) := \mathbf{x}\xsupp [ \{ y_i \ | \ i \in \omega \} ]$.

\item Define $FP_k( \mathbf{x} ) := FP( \mathbf{x'} )$ where $x_n' = x_{n+k}$ for all $n$.


\item $FP$-sets have a natural partial subsemigroup structure induced by $\matf$, i.e., $( \prod_{i\in s} x_i) \cdot (\prod_{i\in t} x_i)$ is defined as in $S$ but only if $\max(s) < \min(t)$. With respect to this restricted operation define $FP^\infty ( \mathbf{x } ) := \delta FP( \mathbf{x} ) = \bigcap_{k \in \omega} \bbar{ FP_k(\mathbf{x}) }$.

\item If the semigroup is written additively, we write $FS( \mathbf{x } )$ etc. accordingly (for finite sums); for $\matf$ we write $FU( \mathbf{x } )$ etc. (for finite unions).
\end{itemize}
\end{Definition}
Instead of saying that a sequence has certain properties it is often convenient to say that the generated $FP$-set does.

The following classical result is the starting point for most applications of algebra in the Stone-\v Cech compactification. We formulate it for partial semigroups.

\begin{Theorem}[Galvin-Glazer Theorem]\label{thm:galvin-glazer} \index{Galvin-Glazer Theorem}\index{Theorem!Galvin-Glazer}
Let $(S,\cdot)$ be a partial semigroup, $p\in \delta S$ idempotent and $A \in p$.
Then there exists $ \mathbf{ x } =  ( x_i )_{i\in \omega}$ in $A$ such that
\[
FP( \mathbf{x } ) \subseteq A.
\]
\end{Theorem}
\begin{proof} 
This can be proved essentially just like the the original theorem, cf.~\cite[Theorem 5.8]{HindmanStrauss}, using the fact that $\sigma(S) \subseteq p $ to guarantee all products are defined.
 \end{proof} 

An immediate corollary is, of course, the following classical theorem, originally proved combinatorially for $\matn$ in \cite{Hindman74}.

\begin{Theorem}[Hindman's Theorem] \label{thm:hindman} \index{Hindman's Theorem}\index{Theorem!Hindman} \index{Finite Sums Theorem|see{Hindman's Theorem}} \index{Theorem!Finite Sums|see{Theorem!Hindman}}
Let $S = A_0 \cup A_1$. Then there exists $i\in \{0,1\}$ and a sequence $ \mathbf{x} $ such that $FP( \mathbf{x } ) \subseteq A_i$.
\end{Theorem}

\section{Union Ultrafilters}

This paper deals primarily with ultrafilters on the partial semigroup $\matf$. The following definitions enable us to speak about the relevant properties of condensations in $\matf$ smoothly.

\begin{Definition}[Condensation, ordered, meshed] \Marginnote{\normalsize Condensation, ordered, meshed} \label{def:condensation,support}
Let $s,t \in \mathbf{F}$ and $\mathbf{s}=(s_i)_{i<N}$  be a disjoint sequence in $\matf $ with $N \leq \omega $.
\begin{itemize}
\item We say that the pair $(s,t)$ is \emph{ordered}, in short $s < t$, if $\max(s) < \min(t) $. 
\item $\mathbf{s}$ is called \emph{ordered}\index{$\matf$!$<$}\index{$\matf$!ordered} if
$ s_i < s_j$ for all $i<j < N$.

\item $v,w \in \matf$ are said to \emph{mesh}, in short $v \sqcap w$, if neither $v<w$ nor $w<v$. \index{$\sqcap$}\index{$\matf$!mesh}
\end{itemize}
\end{Definition}

The following three kinds of ultrafilters were first described in \cite{Blass87-1}.

\begin{Definition}[union ultrafilters] \Marginnote{ \normalsize Union ultrafilters}
An ultrafilter $u$ on $\matf$ is called
\begin{itemize}\index{ultrafilter!union}\index{ultrafilter!{ordered union}}\index{ultrafilter!{stable union}}
\item \emph{union} if it has a base of $FU$-sets.
\item \emph{ordered union} if it has a base of $FU$-sets from ordered sequences.
\item \emph{stable union} if it is union and for any sequence $(FU(\mathbf{s}^\alpha))_{\alpha < \omega}$ in $u$ there exists $FU(\mathbf{t}) \in u$ such that
\[
(\forall \alpha < \omega)\ \mathbf{t} \sqsubseteq^* \mathbf{s}^\alpha,
\]
i.e., $\mathbf{t}$ almost condenses all the sequences $\mathbf{s}^\alpha$ at once.
\end{itemize}
\end{Definition}

It is obvious yet important to note that $FU$-sets always have unique products and all products are defined. At this point it is also good to check the following. Union ultrafilters are elements of $\delta \matf$ and it is not difficult to check that they are idempotent, cf.~\cite[Proposition 3.3]{Blass87-1}, \cite[Theorem 12.19]{HindmanStrauss}. Even though the partial operation on $\matf$ was defined only for disjoint elements it could just as well have been restricted to ordered elements. Of course, this would significantly change the operation on $\matf$ (for example it would not be commutative anymore), but $\sigma(\matf)$ would still be the same and therefore $\delta \matf$. Additionally, the operation on $\delta \matf$ is not changed -- it is after all still the extension of $\cup$ (or $\Delta$) to $\beta \matf$.

The following notion was introduced in \cite{BlassHindman87} to help differentiate union ultrafilters.

\begin{Definition}[Additive isomorphism]  \Marginnote{ \normalsize Additive isomorphism} \label{def:add. isomorphism}
Given partial semigroups $S,T$, call two ultrafilters $p \in \beta S, q\in \beta T$ \emph{additively isomorphic} \index{additive isomorphism} if there exist $FP( \mathbf{ x }) \in p, FP(\mathbf{ y} )\in q$ both with unique products such that the following map maps $p$ to $q$
\[
\phii : FP( \mathbf{ x }) \rightarrow  FP(\mathbf{ y }), \prod_{i \in s} x_i \mapsto \prod_{i  \in s } y_i .
\]
We call this map \emph{the natural (partial semigroup) isomorphism} for $FP$-sets. As mentioned, it extends to a homomorphism (in fact, isomorphism) between $FP^\infty( \mathbf{s} ) $ and $FP^\infty( \mathbf{t} ) $.
\end{Definition}

This notion is a special case of equivalence in the Rudin-Keisler order, but arguably the natural notion for union ultrafilters since every idempotent ultrafilter is isomorphic to an ultrafilter that is not idempotent. For an example, consider the map $\mathbb{F} \rightarrow \mathbb{F}, s\mapsto \{ \max(s) \}$; its extension to $\delta \mathbb{F}$ does have a product of ultrafilters in its range since the set of singletons does not contain any non-trivial products, i.e., union of two disjoint elements.

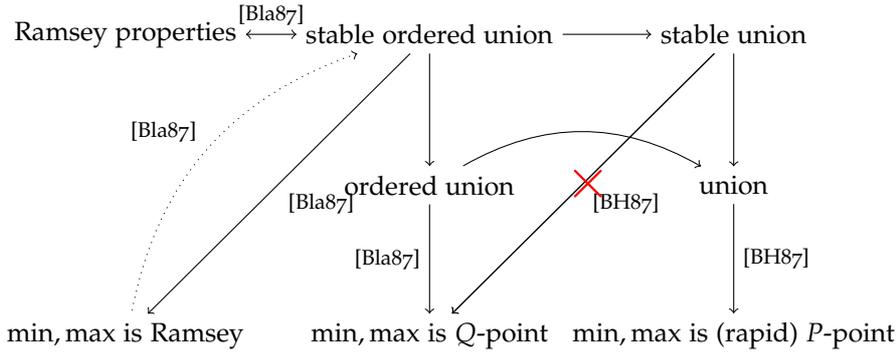
\begin{figure}[h] \caption{Union ultrafilters}\label{fig:union ufs}
\begin{tikzpicture}[scale=2] \footnotesize
[inner sep=2mm]
\node[] (sou) at (0,0) {stable ordered union};
\node[]  (su) at (2,0)  {stable union};
\node[]  (rp) at (-2,0) {Ramsey properties};
\node[]  (ou) at (0,-1) {ordered union};
\node[]  (u)  at (2,-1) {union};
\node[]  (r)  at (-2,-2) {$\min, \max$ is Ramsey};
\node[]  (q)  at (0,-2) {$\min, \max$ is $Q$-point};
\node[]  (p)  at (2,-2) {$\min, \max$ is (rapid) $P$-point};

\path[<->]  (rp) edge node[auto] {\scriptsize \cite{Blass87-1}}    (sou);
\path[->]  (sou) edge node[auto] {}    (ou);
\path[->]  (sou) edge node[auto] {}    (su);
\path[->]  (su) edge node[auto] {}    (u);
\path[->]  (ou) edge[bend left] node[auto] {}    (u);
\path[->]  (ou) edge node[auto,swap] {\scriptsize \cite{Blass87-1}}    (q);
\path[->]  (sou) edge[auto] node[auto] {\scriptsize \cite{Blass87-1}}    (r);
\path[dotted,arrows=<-]  (sou) edge[bend right] node[auto,swap] {\scriptsize \cite{Blass87-1}}    (r);
\path[->]  (u) edge node[auto] {\scriptsize \cite{BlassHindman87}}    (p);
\path[->] (su) edge[] node[auto] {\scriptsize \cite{BlassHindman87}}    (q);
\path[->] (su) edge[] node[] {\LARGE \color{red} \,$\mathbf{\times}$}    (q);

\end{tikzpicture}
\end{figure}

Figure \ref{fig:union ufs} recapitulates the known implications between the types of union ultrafilters with references; the Ramsey properties will be described later. The dotted arrow represents the following: under $CH$, given two non-isomorphic Ramsey ultrafilters, there exists a stable ordered union ultrafilter that maps to them via $\min$ and $\max$.

The one interesting non-implication missing is that a fortiori there consistently exist union ultrafilters that are not ordered union. However, the construction in \cite{BlassHindman87} does not give any direct information on what it means to be an unordered union ultrafilter. In a manner of speaking it is a sledge hammer smashing orderedness so badly that is difficult to identify how orderedness actually fails.
Because of this the main result is dedicated to understanding unordered union ultrafilters. In particular, our construction answers a question by Andreas Blass if there can be unordered union ultrafilters that map to Ramsey ultrafilters via $\max$ and $\min$.
The following result due to Andreas Blass will be needed later.

\begin{Theorem}[Homogeneity (Blass)]\label{thm:blass homogeneity} \index{ultrafilter!selective}\index{ultrafilter!Ramsey}\index{$f$-alternating}
Let $p_0, p_1$ be non-isomorphic, selective ultrafilters and let $f\in 2^\omega$.

If $\mathp(\omega)$ is partitioned into an analytic and coanalytic part, there are $X_i \in p_i$ ($i=0,1$) such that the set of (ranges of) increasing sequences
\[
(x_n)_{n\in \omega} \mbox{ with } x_n \in X_{f(n)}
\]
is homogeneous. We call such sequences \emph{$f$-alternating}.
\end{Theorem}
\begin{proof} 
This is \cite[Theorem 7]{Blass88}
\end{proof}

Regarding this theorem, the following folklore observation is very useful later; cf.~\cite[Lemma 1.2]{Blass87-1}.

\begin{Remark}\label{rem:alternating sets in selectives} \index{ultrafilter!selective}\index{ultrafilter!Ramsey}\index{$f$-alternating}
Given any $f\in 2^\omega$ and non-isomorphic, selective ultrafilters $p_0, p_1$ and $A_i\in p_i$ ($i \in2$), there exists an $f$-alternating sequence in $A_0 \cup A_1$ such that its alternating parts are sets in $p_0$ and $p_1$ respectively.partition $\omega$ in intervals as follows: 
\end{Remark}
\begin{spoiler}
We modify the argument from \cite[Lemma 1.2]{Blass87-1} for alternating sequences using a standard argument for not nearly coherent filters.
\end{spoiler}

\begin{proof}
\plist
\item Pick $A_0 \in p_0$, $A_1 \in p_1$.
\item Let us say $f$ \emph{switches at $i$} if $f(i-1) = j $ and $f(i)=1-j$ (for $j\in 2$). 
\item It is easy to inductively partition $\omega$ into intervals $I_{n}$ large enough so that both $|I_n \cap A_0|$ and $|I_n \cap A_1|$ are at least as large as the longest constant sequence in the range of $f$ up to the $(2n)$th switch; in other words, we make the intervals long enough so that when we can build an $f$-alternating sequence with each alternating ``block'' contained in one of the $I_n$.
\begin{prose}
We will now thin out the ultrafilter sets so that they alternate (though not yet $f$-alternate) but with a twist: the thinned out sets will never meet the same interval $I_n$. After we accomplish this we can fill elements back in from the original $A_i$ to get an $f$-alternating sequence. Since this only enlarges our sets, we stay in our filters.
\end{prose}
\item Consider the interval collapsing map, mapping elements in $I_n$ to $n$.
\item Since this collapsing map cannot map our non-isomorphic selectives to the same ultrafilter, we can find $B_0\in p_0, B_1 \in p_1$ with $B_i \subseteq A_i$ ($i\in 2$) such that $B_0, B_1$ never meet the same interval $I_n$.
\item Next, consider the map on $B_0$ defined by taking $x \in B_0$ to the largest $y\in B_1$ with $y < x$; if this fails map $x$ to $0$.
\item Since $p_0$ is selective, this map becomes injective on a set $C_0\in p_0$, $C_0 \subseteq B_0$. As a result, there must be at least one element from $B_1$ between every two elements in $C_0$.
\item The same argument for $B_1$ (comparing it to $C_0$) yields $C_1\in p_1$, $C_1 \subseteq B_1$ such that there is at least one element from $C_0$ between every two elements in $C_1$.
\item This might, of course, have ruined said property of $C_0$, but we can safely fill in extra elements from $B_1$ to $C_1$ to reestablish it; we still call that possibly larger set $C_1$. In other words, $C_0$ and $C_1$ alternate.
\item Of course, this also means the intervals $I_n$ that are met by $C_0$ and $C_1$ alternate.
\item Finally, by choice of the $I_n$, we can now form an $f$-alternating sequence that incluces $C_0 \cup C_1$.
\plist
\item Simply add elements from the original $A_i$ to the $C_i$ ($i \in 2$) in such a way that they still meet the same intervals but in a block large enough to become $f$-alternating.
\item Since $C_0 \in p_0$, $C_1 \in p_1$, this $f$-alternating sequence is as desired.
\pliste
\pliste
\end{proof}

\section{Stability}

Andreas Blass laid the foundation for all further research regarding union and hence strongly summable ultrafilters in \cite{Blass87-1}. The final theorem from that paper gives a list of potent characterizations of the strongest notion, stable ordered union ultrafilters. (Un)fortunately, not all union ultrafilters are ordered. The first example was constructed in \cite{BlassHindman87} and we will construct an example in the second part. However, all known constructions of union ultrafilter yield stable ones.

In this section we will discuss which of the characterizations for stable ordered union ultrafilters also hold for stable union ultrafilters. 
Because this requires a few definitions that are not relevant for the second half of this paper, we will proceed as follows. We will introduce the one notion that is also of interest for the second part and continue to prove the main result of this section. Following the proof we will discuss the other notions less formally since this does not require as much proof.

\subsection{Stability and the Ramsey property for pairs}

\begin{Definition}[Ramsey property for pairs] \Marginnote{ \normalsize Ramsey property for pairs} Consider $u \in \delta \matf$.
\begin{itemize}

\item We denote the ordered ordered pairs by $\matf^2_<$, i.e.,
\[
\matf^2_< := \{ (s,t)  \in \matf^2 \ |\ s < t \}.
\]
Often $(s<t)$ is a convenient notation for elements in $\mathbb{F}^2_<$.

\item $u$ has the \emph{Ramsey property for pairs}\index{Ramsey property for pairs} if for any finite partition of $\matf^2_<$ there exists $A\in u$ such that $A^2_<$ is homogeneous.

\end{itemize}
\end{Definition}

In \cite[Theorem 4.2]{Blass87-1} Andreas Blass showed that for ordered union ultrafilters the Ramsey property for pairs (and other properties we discuss later) is equivalent to stability. The following result shows that orderedness is not necessary for this equivalence. 

However, it must be stressed that even though the formulation of the Ramsey property is the same, the result is quite different for the unordered case. For an ordered union ultrafilter we get homogeneity for all pairs from the generating sequence. In the unordered case, we do not get such a full property as we cannot check pairs of generators that mesh. We might try to blame this on our formulation of the Ramsey property. Why not ask for partitions of disjoint ordered pairs instead? Unfortunately, this is not possible as the partition of the disjoint pairs into ordered and unordered pairs yields a counterexample for all union ultrafilters, in fact, all idempotent ultrafilter in $\delta \matf$. Every $FU$-set yields both ordered and unordered pairs no matter how nicely the generating sequence behaves.

\begin{Theorem}\label{thm:stable union}\index{ultrafilter!{stable union}}\index{Ramsey property}
A union ultrafilters is stable if and only if it has the Ramsey property for pairs.
\end{Theorem}

\begin{spoiler}
The argument (necessarily) follows the same strategy as the proof of \cite[Theorem 4.2]{Blass87-1}. The forward direction is similar to the proof of Ramsey's Theorem using a non-principal ultrafilter. To get a homogeneous set actually in the ultrafilter stability and a new kind of parity argument is applied.

The reverse conclusion is just as in the original proof by Andreas Blass.
\end{spoiler}

\begin{proof} \plist
\item The Ramsey property for pairs implies stability.
\plist
\item Given any sequence $(FU(\mathbf{s}^\alpha))_{\alpha<\omega}$ in $u$ consider the following set of ordered pairs
\[
\{ (v,w) \in \matf^2_< \ | \  w \in \bigcap_{\alpha < \max(v)} FU(\mathbf{s}^\alpha) \}.
\]

\item Any $FU(\mathbf{t}) \in u$ will yield ordered pairs that are in the above set.
\plist
\item  Pick any $v \in FU(\mathbf{t})$.

\item Take $w>v$ from $FU(\mathbf{t}) \cap \bigcap_{\alpha < \max(v)} FU( \mathbf{s}^\alpha ) \in u$.
\pliste

\item Therefore, by the Ramsey property for pairs, there must be a set $FU(\mathbf{s}) \in u$ such that all ordered pairs are included in the above set.

\item Then $\mathbf{s} \sqsubseteq^* \mathbf{s}^\alpha$ for all $\alpha < \omega$.
\plist
\item Given $\alpha <\omega$, pick $s_i$ with $\max(s_i) > \alpha$.

\item Then all but finitely many $s_j$ have $s_j > s_i$.

\item For such $s_j$ of course $(s_i,s_j) \in FU(\mathbf{s})^2_<$, hence
\[
s_j \in \bigcap_{\beta < \max(s_i)} FU(\mathbf{s}^\beta),
\]
\item  In particular $s_j \in FU(\mathbf{s}_\alpha)$ -- as desired.
\pliste
\pliste 

\item Stability implies the Ramsey property for pairs.

\plist
\item Assume that $A_0 \dotcup A_1 = \matf^2_<$.

\item \Marginnote{Always pick one colour beyond $x$}
Since $u$ is an ultrafilter (in $\delta \matf$)
\[
(\forall x\in \matf)(\exists i) \{ y\in \matf \ |\ (x<y) \in A_i \} \in u.
\]

\item \Marginnote{$A, C_x$ -- almost always pick the same colour.}
Since $u$ is an ultrafilter it concentrates on one colour; without loss it is $0$, i.e., there is $\hl{A\in u}$ such that
\[
(\forall x\in A)\ \hl{C_x := \{ y \in \matf \ |\ (x<y) \in A_0 \}} \in u.
\]

\item Since $u$ is union there are $FU(\mathbf{s}^\alpha) \in u$ (for $\alpha < \omega$) such that
\[
\hl{FU(\mathbf{s}^\alpha)} \subseteq \bigcap_{\max(x) \leq \alpha} C_x.
\]

\begin{prose}
This small simplification ensures that for $x\in A$ we get $FU(\mathbf{s}^{\max(x)}) \subseteq C_x$ by choice of $\mathbf{s}^\alpha$.
\end{prose}

\item \Marginnote{Stability -- almost always pick from the same set}
Since $u$ is stable by assumption, there is $ \hl{FU(\mathbf{s})} \in u$ such that
\[
	{\mathbf{s} \sqsubseteq^* \mathbf{s}^\alpha } \mbox{ for all } \alpha < \omega.
\]

\begin{prose}
Next we introduce a function $j$ that essentially just checks how many members of $\mathbf{s}$ are not included in $FU(\mathbf{s}^\alpha)$.
\end{prose}

\item \Marginnote{$j$ -- Counting where $\mathbf{s}$ fails}
Consider the following function
\[
\hl{j: \omega \rightarrow \omega}, \alpha \mapsto \max \{ \max(s_i)\ |\ s_i \notin FU(\mathbf{s}^\alpha)\}.
\]
Without loss, $j$ is strictly increasing.
\begin{prose}
We can make $j$ strictly increasing by replacing $FU(\mathbf{s}^\alpha) $ with $\bigcap_{\beta \leq \alpha} FU(\mathbf{s}^\beta)$ in the definition of $j$. Alternatively, intersecting $FU(\mathbf{s}^{\alpha+1}) $ with $FU(\mathbf{s}^{\alpha}) $ when we defined them also guarantees $\mathbf{s}^{\alpha+1} \sqsubseteq \mathbf{s}^\alpha $ for all $\alpha < \omega$. In either case, the ``losses'' will at most increase with increasing $\alpha$.
\end{prose}

\item Observe that for all $x\in FU(\mathbf{s})$
\[
\min(x) > j(\alpha) \Rightarrow x \in FU(\mathbf{s}^\alpha).
\]
\plist
\item For $s_i$ this follows by contraposition from the definition of $j$. 
\item Therefore if $\min(x) > j(\alpha)$ this argument implies that all $s_i \subseteq x$ are in $FU(\mathbf{s}^\alpha)$.
\item In particular, so is their union, i.e., $x$.
\pliste

\begin{prose}
After this observation the next goal is to construct $A' \in u$ for which $v<w $ in $A'$ implies $\min(w) > j(\max(v))$. For then $w \in FU(\mathbf{s}^{\max(v)}) \subseteq C_v$. For this a new partition argument is needed.
\end{prose}

\item \label{thin1} \Marginnote{Thinning out 1 -- bounding $j$}
$\{ x \in FU(\mathbf{s}) \ |\ j(\min(x)) < \max(x)\} \in u$.
\plist
\item In any condensation of $\mathbf{s}$ there are $x, x'$ and $x'\cup x$ such that $x < x'$  and $j(\min(x)) <  \max(x')$. 
\item But then calculate
\[
j(\min(x\cup x')) = j(\min(x)) < \max(x') = \max(x\cup x ').
\]
\item Hence any set in $u$ will intersect the above set; so it lies in $u$.
\pliste
\item In particular, there exists $\hl{FU(\mathbf{t}) \in u}$ included in the above set.

\item \Marginnote{Thinning out 2 -- splitting points}
For $x \in FU(\mathbf{t})$ say that \emph{$x$ splits at $n\in x$}\index{splitting points}, whenever
\begin{align*}
& x\cap (n+1), x \setminus (n+1)  \in FU(\mathbf{t})  \mbox{ and}\\
& (\exists t_k)\ x\cap (n+1) < t_k < x \setminus (n+1).
\end{align*}

\begin{figure}[h]
\begin{tikzpicture}[scale=.5]
   \coordinate (y) at (0,0);
   \coordinate (x) at (10,0);
    \draw[thick,->] (y) node[left] {$\ldots$} -- (0,0) --  (x) node[right]
    {$\omega$};
\draw[dotted,red] (.5 cm , -15pt) node[left]{$x$} -- (2.5cm,-15pt) node[anchor=north] {$t_{i}$};
\draw[dotted] (2.75 cm , -10pt) -- (3.75 cm,-10pt) node[anchor=north] {$t_{k}$};
\draw[dotted,red] (4 cm , -15pt) -- (6.5cm,-15pt) node[anchor=north] {$t_{j}$};
\draw[dotted,red] (6 cm , -10pt) -- (8cm,-10pt) node[anchor=north] {$t_{l}$};
\end{tikzpicture}
\caption{Splitting point -- an example}
\end{figure}
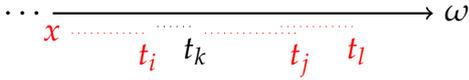

Let $\pi(x)$ be the number of splitting points of $x$, i.e.,
\[
\pi(x) := |\{ n \in x \ |\ x \mbox{ splits at } n\}|
\]

\begin{prose}
The splitting points tell how often $x$ splits into two ordered parts (the one up to $n$ and the one beyond $n$) -- but more importantly with a gap in between.
\end{prose}

\item $\{ x \in FU(\mathbf{t})\ |\ \pi(x)= 1 \bmod 2  \} \in u$.
\plist
\item Any condensation of $t$ will contain some $x < y < z$ and $x \cup z$.
\item In that case, the number of splitting points of $x \cup z$ is
\[
\pi(x\cup z) = \pi(x) + \pi(z) +1.
\]
\item In particular, the number of splitting points for at least one of $x, z, x\cup z$ must be odd.
\pliste
\item In particular, there exists $FU(\mathbf{v}) \in u$ contained in the above set.

\item For any $w_0 < w_1$ in $FU(\mathbf{v})$, there exists $t_j$ with
\[
w_0 < t_j < w_1
\]
\plist
\item Or else $\pi(w_0 \cup w_1) = \pi(w_0) + \pi(w_1)$ would be an even number of splitting points.
\pliste

\item $FU(\mathbf{v})$ is homogeneous, i.e., $FU(\mathbf{v})^2_< \subseteq A_0$.\Marginnote{The conclusion}
\plist
\item For this pick any $w,w' \in FU(\mathbf{v})$ with $w < w'$.

\item By the last step, there exists some $t_j$ with $w < t_j < w'$. Therefore
\[
\min(w') > \max(t_j) > j(\min(t_j)) > j(\max(w))
\]
where the third inequality holds because of step \ref{thin1}, the last because $j$ is strictly increasing.

\item But as we noted just before step \ref{thin1}, this implies $w'\in C_w$, i.e., $(w,w') \in A_0$ -- as desired.
\pliste
\pliste
This concludes the proof.
\pliste \end{proof}

\subsection{Stability and other partition properties}

Let us now discuss the other properties from \cite[Theorem 4.2]{Blass87-1}.

\paragraph{The Ramsey property for $k$-tuples} It is straightforward to generalise the Ramsey property for pairs to $k$-tuples (with $k < \omega$) as follows. An ultrafilter $u$ has the \emph{Ramsey property for $k$-tuples} if for every partition of $\mathbb{F}^k_< := \{ (s_0,\ldots, s_{k-1}) \ |\ (\forall i<k-1)s_i < s_{i+1} \}$ we can find $A\in u$ such that $A^k_<$ is homogeneous. It is not difficult to show by an induction much like the induction used for Ramsey's Theorem for $\omega$ that the Ramsey property for pairs implies the Ramsey property for $k$-tuples for all $k < \omega$. An alternative argument follows from the property described in the next paragraph.

\paragraph{The infinitary Ramsey property} In \cite[Theorem 4.2]{Blass87-1} Andreas Blass also discusses the infinitary analogue of the Ramsey properties. For this consider the set of ordered $\omega$-sequences $\mathbb{F}^\omega_< := \{ \mathbf{s} \in \mathbb{F}^\omega \ | \ (\forall i < \omega) s_i < s_{i+1} \}$. Then $u$ has \emph{the infinitary Ramsey property} if for every partition of $\mathbb{F}^\omega_<$ into an analytic and co-analytic part there exists $A\in u$ such that all ordered subsets of $A$ are in the same part. It is not difficult to check that the proof in \cite{Blass87-1} does not require the union ultrafilter to be ordered. It might be worthwhile to check that the strength of this infinitary partition property suffers even more than the finitary ones from dropping the ordered union requirement. For a stable ordered union ultrafilter not only do we get the infinitary Ramsey property, but the homogeneous set itself is generated by an ordered sequence, hence that ordered sequence is in that part of the partition. In the unordered case this statement simply does not make sense as the partition only covers ordered sequences.

\paragraph{Characterization via $\min$ }
The last two properties from \cite[Theorem 4.2]{Blass87-1} are formulated in terms of ultrapowers of $\omega$. To keep our discussion short, we assume some basic knowledge of ultrapowers of $\omega$; for a concise introduction cf.~\cite[Section 1]{Blass87-1} (which is available at \cite{BlassWWW}). Given an ultrafilter $V$ on a countable set $I$, we say that $f,g \in I^\omega$ are in the same \emph{sky} if there exists $h,h' \in \omega^\omega$ such that on a set in $V$ we have $g \leq h \circ f $ and $f \leq h'\circ g $. It is known that equivalently $f,g$ are in the same sky if there exist finite-to-one $h,h'$ such that $h\circ f = h' \circ g$. Skies are obviously order convex.

For lack of a better term, we say that $u$ is \emph{stable via $\min$} if whenever $f \in \mathbb{F}^\omega, g\in \omega^\omega$ and $A\in u$ such that $f(s) < g\circ \min(s)$ for all $s\in A$, then there exists $h \in \omega^\omega$ such that $f(s) = h\circ \min(s)$ on some set in $u$; in other words, then the values of $f(s)$ only depend on $\min(s)$.\footnote{In terms of the ultrapower this means that $\min$ generates an initial segment of the ultrapower.} 

In \cite[Theorem 4.2]{Blass87-1} it is 
essentially\footnote{''Essentially'' in the sense that the relevant part of \cite[Theorem 4.2]{Blass87-1} includes the condition that the image of the union ultrafilter under $\min$ is a P-point. This fact was established later for all union ultrafilters in \cite[Theorem 2]{BlassHindman87}.} 
shown that this notion is equivalent to stability for ordered union ultrafilters. For union ultrafilters we can show two things. On the one hand, the following observation shows that stability via $\min$ implies stability. On the other, the next section will include an example showing that (consistently) stability via $\min$ does not imply that a union ultrafilter is ordered.

\begin{Theorem}
If a union ultrafilter is stable via $\min$ then it has the Ramsey property for pairs. In particular, it is stable.
\end{Theorem}
\begin{spoiler}
We only sketch the argument since it is just a recombination of the argument for ordered union ultrafilters with a recent result by Andreas Blass.
\end{spoiler}
\begin{proof}
\plist
\item Let $u$ be a union ultrafilter that is stable via $\min$.
\item The Ramsey property for pairs is easily seen to be equivalent to the statement that $\mathbb{F}^2_<$ and $(A\times A)_{A \in u}$ generate an ultrafilter on $\mathbb{F}^2$, namely the tensor product $u\otimes u$.
\item So take any ultrafilter $V$ on $\matf^2$ containing all these sets. We show that $V= u \otimes u$.
\item By a characterization of tensor products due to Puritz \cite{Puritz1972215}, cf.~also \cite[Section 1]{Blass87-1}, it suffices to show that whenever $g_1,g_2 \in \omega^\omega$, then $g_1 \circ \pi_1$ lies in a lower sky than $g_2 \circ \pi_2$ (unless the latter is constant on a set in $V$).\footnote{Where $\pi_i$ is the projection to the $i$-th coordinate ($i\in 2$).}
\item It suffices to compare $\max \circ \pi_1$ with $\min \circ \pi_2$.
\plist 
\item $\pi_2(V)=u$ and by \cite[Theorem 2]{BlassHindman87} $\min(u)$ is a P-point, i.e., the sky of $\min$ contains exactly two skies, one of them the sky of constant functions. 
\item Combining this with stability via $\min$, we get that $\min \circ \pi_2$ is in the lowest non-standard sky for elements of the form $g_2 \circ \pi_2$. 
\item Also, $\max$ is finite-to-one, so $\max \circ \pi_1$ is in the highest sky for elements of the form $g_1 \circ \pi_1$.
\pliste
\item $\max \circ \pi_1$ is at most in the same sky as $\min \circ \pi_2$.
\plist
\item By assumption on $V$ we have $\max \circ \pi_1 < \min \circ \pi_2$ on $\matf_<^2 \in V$.
\pliste
\item $\max \circ \pi_1$ is not in the same sky as $\min  \circ \pi_2$.
\plist
\item By \cite[Theorem 38]{Blass20092581}, $\min(u), \max(u)$ are not near coherent filters, i.e., no two finite-to-one maps will map $\min(u)$ and $\max(u)$ to the same ultrafilter.
\item But by \cite[Theorem 2]{BlassHindman87} $\min(u), \max(u)$ are both $P$-points.
\item So no two maps will map $\min(u)$ and $\max(u)$ to the same non-principal ultrafilter.
\item In other words, on any set in $V$, we have $h \circ \max \circ \pi_1 \neq h' \circ \min \circ \pi_2$ for any $h,h' \in \omega^\omega$ (unless both sides are constant).
\item So $\max \circ \pi_1$ and $\min  \circ \pi_2$ are not in the same sky.
\pliste 
\item Therefore $V = u\otimes u$, i.e., $u$ is stable.
\pliste
\end{proof}

\paragraph{The canonical partition property} We say that that an ultrafilter $u\in \delta \matf$ has \emph{the canonical partition property} if for each $f:\mathbb{F} \rightarrow \omega $ there exists $A\in u$ such that $f\restr A$ has one of the following properties: 
\begin{itemize}
\item $f\restr A$ is constant,
\item $f\restr A = g\circ \min \restr A$ for some injective $g \in \omega^\omega$, in particular, the values only depend on $\min(s)$,
\item $f\restr A = g \circ (\min,\max) \restr A$ for some injective $g:\omega^2 \rightarrow \omega$, in particular, the values of $f\restr A$ depend only on $(\min(s), \max(s))$,
\item $f\restr A = g \circ \max \restr A$ for some injective $g\in \omega^\omega$, in particular, the values of $f\restr A$ depend only on $\max(s)$.
\item $f\restr A$ is injective.
\end{itemize}
Again, the proof of \cite[Theorem 4.2]{Blass87-1} not only shows that the canonical partition property implies stability via $\min$, hence stability, for ordered union ultrafilters, but also for union ultrafilters in general. It remains open whether this property is equivalent to stability. Also, we do not know if it implies orderedness.

\subsection{Stability and additive isomorphisms}

We end this section with the following application of stability which will be useful later. The result is remininiscent of the role of $P$-points and Ramsey ultrafilters in the Rudin-Keisler order.

\begin{Lemma}[Stability and Additive Isomorphisms] \label{lem:stability vs add. homomorphisms}
Every additively isomorphic image of a stable union ultrafilter is a stable union ultrafilter and every additively isomorphic image of a stable ordered union ultrafilter is a stable ordered union ultrafilter.

\end{Lemma}

\begin{spoiler}
Stability is straightforward; for orderedness we use the Ramsey property of stable ordered union ultrafilters.
\end{spoiler}

\begin{proof} \plist
\item Let $u'$ and $u$ be additively isomorphic ultrafilters, i.e., there exist $FU(\mathbf{s}) \in u$, $FU(\mathbf{x}) \in u'$ such that
\[
\pi: FU(\mathbf{s}) \rightarrow FU(\mathbf{x}), \prod_{i\in F} s_i \mapsto \prod_{i\in F} x_i
\]
additionally has $\pi(u)=u'$.

\item \label{stable} If $u$ is stable, so is $u'$.
\plist
\item Given a sequence of condensations with $FU$-sets in $u'$ we may assume without loss that all of them condense $\mathbf{x}$.

\item But then the preimages under $\pi$ form a sequence of condensations in $u$. 

\item Applying the stability of $u$ yields a common almost condensation of the images.

\item Then its image under $\pi$ is exactly the desired common almost condensation in $u'$.
\pliste

\item If $u$ is stable ordered, so is $u'$.
\plist
\item By step \ref{stable}, we only need to show that $u'$ is ordered. Pick any $A\in u'$; we may assume without loss $A\subseteq FU(\mathbf{x})$ and also that $\mathbf{s}$ (from the definition of $\pi$) is ordered.

\item \Marginnote{Working in $u$: a partition for orderedness.}
Consider
\[
X:= \{ (v, w) \in [\pi^{-1}[A]]_<^2 \ | \ \max(\pi(v)) < \min(\pi(w))\}.
\]

\item By the Ramsey property from Theorem \ref{thm:stable union}, there exists an ordered sequence $t$ such that $FU(\mathbf{t}) \in u$ and $FU(\mathbf{t})$ is either included in or disjoint from $X$.

\item But it cannot be disjoint from $X$.
\plist
\item Since $\pi$ is injective, there exist $t_i < t_j$ such that $\pi(t_i) < \pi(t_j)$.
\item In fact, given any $t_i$ all but finitely many $t_j$ have this property.

\item Then $( t_i, t_j) \in X \cap FU(\mathbf{t})^2_<$.
\pliste

\item But this implies that $\pi[FU(\mathbf{t})]=FU(\pi[\mathbf{t}])$ is ordered -- and of course in $u'$ and refining $A$.
\item In other words, $u'$ is ordered union.
\pliste
This concludes the proof.
\pliste \end{proof}

\section{Unordered union ultrafilters}

We now turn to our main result that selectivity of the image under $\min$ and $\max$ cannot indicate orderedness of a union ultrafilter. At first this is a negative result since the alternative would probably have involved a new partition theorem involving Ramsey ultrafilters.
However, the construction of the counterexample offers an answer to the simple question: What does an unordered union ultrafilter look like? As mentioned earlier the construction in \cite{BlassHindman87} does not really answer this question. Nevertheless, its proof represents a blueprint for constructions of (stable) union ultrafilters.

By definition, to be unordered means that there must be a ``special'' $FU$-set in the ultrafilter that will not be refined to an ordered $FU$-set in the ultrafilter.\footnote{We will always have some ordered $FU$-sets in any (union) ultrafilter, e.g., $\matf = FU( (\{ i\} )_{i \in \omega} ) $.} In particular, the sequence generating the $FU$-set itself cannot be ordered. But if a sequence is not ordered, it is meshed in the sense that some of its members must mesh. Of course, such a sequence will be condensed again and again -- and yet no ordered condensation can be allowed. So the question becomes: what might this meshing look like? 

Let us do some handwaving arguments on some simple attempts that are doomed to fail. Since any union ultrafilter is in $\delta \matf$ and our sequence is disjoint, there must be ``arbitrarily late'' meshing, i.e., if only finitely many elements of $\mathbf{s}$ mesh we have already lost. It is also easy to see that union ultrafilters concentrate on condensations that contain unions of many members of the sequence, e.g., because the sequence itself will not be in the union ultrafilter; therefore there cannot be a bound on the number of $s_i$ which mesh. Finally, by parity arguments the meshing cannot be only of, e.g., the form $s_{2i} \sqcap s_{2i+1}$, since a union ultrafilter will concentrate on those with an even number of adjacent indices -- so any union ultrafilter will condense such a sequence to an ordered sequence. Finally, the critical concern will have to be, whether there can be enough meshing while keeping the images under $\min$ and $\max$ Ramsey ultrafilters.

The main result of this section is as follows.

\begin{Theorem}[Stable Unordered  Union Ultrafilters]
Assume \CH. There exists a stable union ultrafilter $u$ with $\min(u)$ and $\max(u)$ selective, but there exists $FU(\mathbf{s}) \in u$ such that for every ordered sequence $\mathbf{t}$
\[
\mathbf{t} \sqsubseteq \mathbf{s} \Rightarrow FU(\mathbf{t}) \notin u.
\]
In fact, any two non-isomorphic selective ultrafilters can be prescribed for $\min$ and $\max$.
\end{Theorem}
Note that the assumption of \CH\ can be weakened to essentially iterated Cohen forcing; this will be discussed at the end of the section.

Fortunately, with the help of the lemma in the previous section this implies a stronger version guaranteeing rigidity under additive isomorphisms.

\begin{Corollary}[Unordered Union Ultrafilters]\index{ultrafilter!{unordered union}}
Assume \CH. There exists a stable union ultrafilter $u$ with $\min(u)$ and $\max(u)$ selective, but $u$ is not additively isomorphic to an ordered union ultrafilter.

In fact, any two non-isomorphic selective ultrafilters can be prescribed for $\min$ and $\max$.
\end{Corollary}
\begin{proof}
This follows from the above theorem and Lemma \ref{lem:stability vs add. homomorphisms}\end{proof}

\subsection{The construction} 

Recall \Marginnote{The critical issue -- a special set $FU(\mathbf{s})$} the goal: however the union ultrafilters $u$ is constructed, it must include a set $FU(\mathbf{s})$, such that \emph{any ordered condensation} $\mathbf{t} \sqsubseteq \mathbf{s} $ is excluded, i.e., $FU(\mathbf{t}) \notin u$.

It is not difficult to put together a union ultrafilter with a base of unordered $FU$-sets. But this does not suffice, since there might be a different base of ordered $FU$-sets by accident.

To prevent this, no unordered condensation that is added in the inductive construction can accidentally be, at the same time, a condensation of some other, ordered condensation of the fixed $FU(\mathbf{s})$ (thus including that ordered condensation of $FU(\mathbf{s})$ in the ultrafilter $u$ as well).

This means that every chosen sequence must eventually have a high degree of meshing not just in itself but due to the $s_i$ that appear in its support. The following definition prepares for the right notion of meshing.

\begin{Definition}[The meshing graph]\label{def:meshing graph} \index{meshing graph}\index{graph!meshing} \Marginnote{\normalsize The meshing graph}
Given $\mathbf{s}=(s_i)_{i<N}$ (for some $N \leq \omega $) and some condensation $\mathbf{t}=(t_j)_{j<K}$ of $\mathbf{s}$ (for some $K\leq N$) we define  the \emph{meshing graph $G_t$} to be the graph on the vertices $\{t_j \ | \ j < K \}$ with edges
\[
E(G_t) = \{ \{t_i, t_j\} \ |\ (\exists s_n \subseteq t_i, s_m \subseteq t_j)\ s_n \sqcap s_m \},
\]
i.e., there is an edge whenever two $t_j$ are meshed and this meshing is caused by two elements from $\mathbf{s}$.
\end{Definition}

This notion allows to discuss the degree of meshing in terms of the connectedness of the graph. On the one hand, it is an advantage to connect to graph theory and graph colourings. On the other hand, it is unclear how well connected the graph should be -- and it is not trivial to get Ramsey-type theorems for graphs that allow a flexible degree of connectedness. Fortunately, it will be enough to work with complete graphs.

To begin the construction, a thoroughly meshed sequence is required. After all, in an inductive construction under \CH, the critical $FU$-set must appear after countably many steps so it might as well appear right away. The general case without a preselected $FU$-set will be discussed later.

\begin{Remark}[Fix the meshed sequence $\mathbf{s}$]\Marginnote{\normalsize The fixed sequence $\mathbf{s}$}
From now on fix a sequence $\mathbf{s}=(s_i)_{i \in \omega}$ such that for any $n$ there exist $i_0<\ldots<i_n$ such that
\[
G_{(s_{i_0}, \ldots, s_{i_n})}
\]
is a complete graph with $n+1$ vertices.\footnote{Here the meshing graph is computed with respect to $\mathbf{s}$ itself.}
\end{Remark}

This simply means that the sequence includes arbitrarily large segments that have the best meshing. As mentioned, for now it is enough to pick any such sequence (which is easy to construct inductively). It will be proved how to find such a sequence with respect to two prescribed selective ultrafilters at the end of the construction. It is useful to note that such a sequence might (and later will) be chosen "nearly ordered" in the sense that the increasingly large complete graphs appear in an ordered fashion.

The following definition tries to capture the right kind of meshing that is needed for condensations and more generally for sets that are suitable for the ultrafilter.

\begin{Definition}[$\mathbf{s}$-meshed] \Marginnote{\normalsize $\mathbf{s}$-meshed} \index{$\mathbf{s}$-meshed}
A set $A\subseteq \matf$ is called \emph{$\mathbf{s}$-meshed} if for any $n \in \omega$ there exist (disjoint) $\mathbf{t}=(t_i)_{i<n}$ such that
\begin{itemize}
\item $FU(\mathbf{t}) \subseteq (FU(\mathbf{s}) \cap A)$
\item The meshing graph $G_t$ is a complete graph.
\end{itemize}
We call such a finite sequence an \emph{$n$-witness} of $A$.\footnote{Note that an $\mathbf{s}$-meshed set is compatible with $\delta \matf$ since for any $v \in \matf$ any disjoint sequence of length $\max(v)+2$ must have an element in $\sigma(v)$.}
\end{Definition}

A set $A$ is $\mathbf{s}$-meshed if there are members of $A$ that have a high degree of meshing and additionally the witnesses for the meshing are given by arbitrarily large, finite $FU$-sets where the members of the $\mathbf{s}$-support mesh very much.

The following observation should support the claim that this is the right notion for this setting, i.e., such sets do not force us to add ordered condensations to an ultrafilter.

\begin{Proposition}
If $A$ is $\mathbf{s}$-meshed, then it is not included in $FU(\mathbf{t})$ for any ordered $\mathbf{t} \sqsubseteq \mathbf{s}$.
\end{Proposition}

\begin{proof} \plist
\item To be an ordered condensation $\mathbf{t} \sqsubseteq \mathbf{s}$ means that $G_t$ has no edges.


\item There are no disjoint elements in $FU(\mathbf{t})$ with a non-empty meshing graph.
\plist
\item Assume $v,w \in FU(\mathbf{t})$ have an edge, i.e., there exist $s_i \subseteq t_j \subseteq v, s_k \subseteq t_l \subseteq w$ with $s_i \sqcap s_j$. 
\item Therefore $t_j \sqcap t_l$.
\item Since $\mathbf{t}$ is ordered, this implies $t_j = t_l$, i.e., $t_j \subseteq v \cap w \notnil $
\pliste
\item Hence $FU(\mathbf{t})$ cannot include an $\mathbf{s}$-meshed set.
\pliste \end{proof}

To be able to link the new notion with ultrafilters it needs to be partition regular. This requires the following classical result which is sometimes called ''finite Hindman's Theorem`` even though it historically preceded and motivated Hindman's Theorem.

\begin{Theorem}[Folkman-Rado-Sanders]\label{thm:finite hindman}
For any $n \in \omega $ there exists $h(n) \in \omega$ such that for any disjoint sequence $\mathbf{x}= (x_i)_{i<h(n)}$ in $\matf$ the following holds:

Whenever \Marginnote{Folkman-Rado-Sanders Theorem} $FU(\mathbf{x})$ is finitely partitioned, there exists a condensation of length $n$ with a homogeneous $FU$-set.
\end{Theorem}

\begin{proof} 
The original discovery is attributed to Folkman and Sanders independently; it follows from Rado's Theorem and from the Graham-Rothschild Parameter Sets Theorem, \cite[cf.~Corollary 3]{GrahamRothschild71}. For a proof from Hindman's Theorem by a compactness argument see \cite[Theorem 5.15]{HindmanStrauss}; for a more recent overview on its combinatorial aspects cf.~\cite{ProemelVoigt90}.
\end{proof}

This allows for the proof of the first piece of the puzzle.

\begin{Lemma}[$\mathbf{s}$-meshed partition regular]\label{lem:meshed partition regular}
The notion of being $\mathbf{s}$-meshed is partition regular.

In particular, any $\mathbf{s}$-meshed set is included in an ultrafilter consisting only of sets that are $\mathbf{s}$-meshed.
\end{Lemma}
\begin{spoiler}
Given a finite partition of an $\mathbf{s}$-meshed set, the Folkman-Rado-Sanders Theorem implies large homogeneous condensations. To get $n$-witnesses it turns out that a homogeneous condensation inherits a complete meshing graph.
\end{spoiler}

\begin{proof} 
\plist
\item Clearly, $\matf$ itself is $\mathbf{s}$-meshed since it contains $FU(\mathbf{s})$. 

\item So fix an arbitrary $\mathbf{s}$-meshed set $A$ and any partition $A= A_0 \dotcup A_1$.

\item Since $A$ contains arbitrarily large witnesses, the Folkman-Rado-Sanders Theorem (plus the pigeon hole principle) implies that in either $A_0$ or $A_1$ there are arbitrarily large condensations of these witnesses.

\item Any condensation $\mathbf{v}$ of a witness $\mathbf{t}$ has a complete meshing graph.
\plist
\item Given $v_i,v_j$ in the condensation, there are $t_k \subseteq v_i, t_l \subseteq v_j$.

\item Since $\mathbf{t}$ has a complete meshing graph, there exists $s_m \subseteq t_k \subseteq v_i$, $s_n \subseteq t_l \subseteq v_j$ with $s_m \sqcap s_n$.

\item Therefore $v_i, v_j$ are connected in the meshing graph $G_\mathbf{v}$.
\pliste
\item Therefore, either $A_0$ or $A_1$ is $\mathbf{s}$-meshed.
\pliste \end{proof}

The next step is to show that the ultrafilters containing $\mathbf{s}$-meshed sets are algebraically rich.

\begin{Lemma}[The meshing semigroup]\label{lem:meshing semigroup}
The set
\[
H := \{ p \in FU^\infty(s)\cap \delta \matf \ |\ (\forall A\in p)\, A \mbox{ is } \mathbf{s}\mbox{-meshed} \}
\]
is a closed subsemigroup of $\delta \matf$.
\end{Lemma}

\begin{proof} \plist
\item $H$ is a closed subset of $\delta \matf$ since it is defined by a constraint on all members of its elements.

\item Lemma \ref{lem:meshed partition regular} implies that it is not empty.

\item $H$ is a subsemigroup.
\plist
\item Pick arbitrary $p,q \in H$ and $V\in p$, $(W_v)_{v\in V}$ in $q$; in particular all these sets are $\mathbf{s}$-meshed.

\item Then $\bigcup_{v\in V} (v \cdot W_v) $ is $\mathbf{s}$-meshed.
\plist
\item Pick any $n \in \omega$.

\item By assumption on $p$ there exists \hl{an $n$-witness $\mathbf{t}=(t_i)_{i<n}$} such that
\[
FU(\mathbf{t}) \subseteq V.
\]
\item Similarly, by assumption on $q$, there exists \hl{an $n$-witness $\mathbf{t'}=(t'_i)_{i<n}$} such that
\[
FU(\mathbf{t'}) \subseteq \bigcap_{x\in FU(\mathbf{t}_0,\ldots,t_n)} W_v \cap \sigma(\bigcup_{i\leq n} t_i) \quad (\in q).
\]

\item But then for \hl{$ \mathbf{v} = (v_i)_{i<n}$} with $v_i := t_i \cup t_i'$ in fact
\[
FU(\mathbf{v}) \subseteq \bigcup_{v\in V} v \cdot W_v.
\]
\item Additionally, $G_v$ is a complete graph since $G_t$ was (or since $G_{t'}$ was) -- making the sets ``fatter'' only increases the chance of being meshed.

\item In particular, the set is $\mathbf{s}$-meshed -- as desired.
\pliste
\item Therefore, $p \cdot q \in H$.
\pliste
This competes the proof.
\pliste \end{proof}

The next step is to show that the preimage filters under $\min$ and $\max$ are compatible with $H$, i.e., contain $\mathbf{s}$-meshed sets.

\begin{Lemma}\label{lem:preimages are meshed}
If $A \cap \min[FU(\mathbf{s})]$, $B \cap \max[FU(\mathbf{s})]$ are both infinite, then
\[
\min\!^{-1}[A] \cap \max\!^{-1}[B]
\]
is $\mathbf{s}$-meshed.
\end{Lemma}

\begin{spoiler}
Pick three sets of members of $\mathbf{s}$: one set to get the prescribed minimum, another set to get the meshing, and finally a set to get the prescribed maximum.
\end{spoiler}

\begin{proof} 
\plist
\item Given $n \in \matn$ we pick three times $n$-many elements of the sequence $\mathbf{s}=(s_i)_{i\in \omega}$.

\item Since $A$ is infinite, it is possible to pick \hl{$(s_{i_k})_{k<n}$ with $\min(s_{i_k}) \in A$}.

\item By the meshing of $\mathbf{s}$, it is possible to pick \hl{$(s_{j_k})_{k<n}$} with a complete meshing graph but lying beyond everything chosen so far.\footnote{In other words, with minima greater than the greatest maximum so far.}

\item Since $B$ is infinite, it is possible to pick \hl{$(s_{l_k})_{k<n}$ with $\max(s_{l_k}) \in B$}, again beyond everything chosen so far.

\item Then \hl{$(t_k)_{k<n}$ defined by $t_k := s_{i_k} \cup s_{j_k} \cup s_{l_k}$ } is an $n$-witness for $\min^{-1}[A] \cap \max^{-1}[B]$.
\pliste \end{proof}

An easy corollary is the following.

\begin{Corollary}
Let $p_1$ and $p_2$ be ultrafilters including $\min[FU(\mathbf{s})]$, $\max[FU(\mathbf{s})]$ respectively. Then
\[
\widebar{\min\, \! ^{-1}(p_1)} \cap \widebar{\max\,\!^{-1}(p_2)} \cap H \notnil.
\]
\end{Corollary}

\begin{proof} 
By Lemma \ref{lem:preimages are meshed} all elements of the preimage filter are $\mathbf{s}$-meshed. A standard application of Zorn's Lemma (or equivalenty compactness) allows us to extend any such filter to an ultrafilter in $H$.
\end{proof}

Note also that $\bbar{ \min\! ^{-1}(p)}$ is a right ideal, $\bbar{ \max \! ^{-1}(p)}$ a left ideal in $\delta \mathbb{F}$ for any $p \in \beta \mathbb{N}$; in particular their intersection is a closed subsemigroup. This is easily checked, cf.~also \cite[Section 2.3]{Krautzberger09} and \cite[Theorem 6.9]{HindmanStrauss}.

\subsection{Main lemma and theorem}

After the preparations are complete it is possible to tackle the main lemma for the inductive construction. Let $\langle \min\!^{-1}(p_1) \cup \max\!^{-1}(p_2) \rangle $ denote the filter generated by the union of the coherent filters $\min\!^{-1}(p_1)$, $\max\!^{-1}(p_2)$.

\begin{Lemma}[Main Lemma]\label{lem:unordered main lemma}\index{ultrafilter!{unordered union}}
Assume we are given non-isomorphic, selective ultrafilters $p_1,  p_2$ with $\max[FU(\mathbf{s})] \in p_1$ and $\min[FU(\mathbf{s})] \in p_2 $ as well as some $X\subseteq \matf$.

For every $\alpha < \omega$ let $\mathbf{t}^\alpha=(t_i^\alpha)_{i\in \omega}$ be a sequence such that
\begin{align*}
& \mathbf{t}^{\alpha +1} \sqsubseteq^* \mathbf{t}^\alpha \\
& FU(\mathbf{t}^\alpha) \mbox{ is $\mathbf{s}$-meshed } \\
& FU(\mathbf{t}^\alpha) \in \langle \min\!^{-1}(p_1) \cup \max\!^{-1}(p_2)\rangle.
\end{align*}

Then there exists $ \mathbf{z} = (z_i)_{i \in \omega} $ such that
\begin{align*}
& \mathbf{z} \sqsubseteq^* \mathbf{t}^\alpha \mbox{ for every } \alpha < \omega, \\
& FU(\mathbf{z}) \subseteq X \mbox{ or } FU(\mathbf{z}) \cap X = \nil, \\
& FU(\mathbf{z}) \mbox{ is $\mathbf{s}$-meshed} \\
& FU(\mathbf{z}) \in \langle \min\!^{-1}(p_1) \cup \max\!^{-1}(p_2) \rangle .
\end{align*}
\end{Lemma}

\begin{spoiler}
By a standard Galvin-Glazer argument there exists a common almost condensation of the given $FU$-sets and $X$ (or its complement). Since all sets are $\mathbf{s}$-meshed, the condensation can be $\mathbf{s}$-meshed. The Homogeneity Theorem \ref{thm:blass homogeneity} ensures such a condensation can be found in $\min\!^{-1}(p_1) \cup \max\!^{-1}(p_2)$.
\end{spoiler}

\begin{proof} 
\plist
\item By the assumptions,
\[
H \cap \widebar{\min\!^{-1}(p_1)} \cap \widebar{ \max\!^{-1}(p_2)} \cap \bigcap_{\alpha < \omega} FU^\infty(\mathbf{t}^\alpha) \notnil.
\]
\item As\Marginnote{$e$ -- the helpful idempotent}  an intersection of closed semigroups it is a closed semigroup which therefore contains an idempotent $\hl{e \in \delta \matf}$. Without loss $X\in e$; in particular $X$ is $\mathbf{s}$-meshed.

\begin{prose}
The aim is to apply the Homogeneity Theorem \ref{thm:blass homogeneity}.
\end{prose}

\item 
Consider the following analytic set in $\mathp(\omega)$.\Marginnote{An analytic set}
\begin{align*}
\{ Y\subseteq \omega \ | \ (\exists \mathbf{z} =(z_i)_{i \in \omega}) &  \ Y = \min[FU(\mathbf{z})] \cup \max[FU(\mathbf{z})] , \\
& (\forall \alpha < \omega )\ z \sqsubseteq^* \mathbf{t}^\alpha, \\
& FU(\mathbf{z}) \subseteq X, \\
& FU(\mathbf{z}) \mbox{ is } FU\mbox{-meshed} \}.
\end{align*}
This set encodes all the candidates for the claim.

\item Define \Marginnote{\small $010011000111\ldots$} \hl{$f \in 2^\omega$} inductively to have $n$-many $0$'s followed by $n$-many $1$'s for each $n$ in increasing order of $n$'s.

\begin{prose}
We will apply Theorem \ref{thm:blass homogeneity} to get an $f$-alternating sequence whose alternating blocks form sets in our selective ultrafilters. An $f$-alternating sequence alternates by picking $n$-many elements from both sets at the $n$-th step. It is useful to check that if we asked for alternating instead of $f$-alternating then we would always miss the analytic set, since any $\mathbf{z}$ with alternating minima and maxima must be ordered, so it isn't $FU$-meshed. On the other hand, a ``nearly ordered'' sequence such as $\mathbf{s}$ comes in ordered blocks of completely meshed finite sequences. The minima and maxima of such a sequence are precisely $f$-alternating which is why we choose $f$-alternating here.
\end{prose}

\item \Marginnote{Using homogeneity}
By Theorem \ref{thm:blass homogeneity} there exist \hl{$Y_1 \in p_1$, $Y_2 \in p_2$} such that
\[
\{ A \subseteq Y_1 \cup Y_2 \ |\ A\ f\mbox{-alternating} \}
\]
is either contained in or disjoint from the analytic set.

\begin{prose}
If the above set is included in the analytic set then Remark \ref{rem:alternating sets in selectives} guarantees the existence of the set desired to complete the proof. Fortunately, given any $Y_1, Y_2$ a Galvin-Glazer argument shows that the set can never be disjoint.
\end{prose}

\item But for every \Marginnote{Applying Galvin-Glazer to find $\mathbf{t}$}
$Y_1\in p_1$, $Y_2 \in p_2$ there exists \hl{$\mathbf{t}=(t_i)_{i \in \omega}$}, a common almost condensation of $(\mathbf{t}^\alpha)_{\alpha<\omega}$ such that
\[
\hl{FU(\mathbf{t}) \subseteq Z } := X \cap \min\!^{-1}[Y_1] \cap \max\!^{-1}[Y_2] \cap FU(\mathbf{s});
\]
additionally, $\mathbf{t}$ is $\mathbf{s}$-meshed and the minima and maxima are $f$-alternating.
\plist
\item First note $Z \in e$; without loss $Z^\star=Z \in e$ (with respect to $e$).

We construct the desired sequence by induction.

\item At the inductive step $n$, having constructed $t_0,\ldots,t_k$ (where $k=\sum_{i=0}^{n-1} i$) we assume by induction hypothesis that the following intersection is in $e$
\[
Z^\star \cap  \bigcap_{x\in FU(\mathbf{t}_0,\ldots,t_k)} x^{-1} Z^\star \cap \sigma(\bigcup_{i<k}t_i) \cap \bigcap_{\alpha<n} FU(\mathbf{t}^\alpha).
\]
\item Pick an $\mathbf{s}$-meshing $n$-witness $t_k,\ldots,t_{n+k}$ from it.\footnote{For later reference, note the following. We have a lot of freedom at this point to impose other properties on these $n+1$-many elements of $\mathbf{t}$. In particular, we can first choose an ordered sequence of length $n+1$, then a sequence of meshing witnesses beyond those such that the union of the $i$th from the ordered part with the $i$th from the meshed part is just as good to continue our induction, i.e., all finite unions are still in $Z$. This kind of ''late meshing'' will be needed for an observation regarding stability via $\min$ at the end of this section. }

\item As usual in the Galvin-Glazer argument, the analogous intersection for $FU(t_0,\ldots. t_{k+n})$ is again in $e$.

\item The resulting sequence is $\mathbf{s}$-meshed by construction.

\item Note that for a sequence of length $n$ with a complete meshing graph, all minima must come before all maxima. Since the witnesses are chosen in an ordered fashion, this implies that the entire sequence has $f$-alternating minima and maxima.
\pliste
\item By Remark \ref{rem:alternating sets in selectives} there exists a sequence $\mathbf{z}$ for $Y_1$ and $Y_2$ themselves --- and with all the desired properties to conclude the proof.
\pliste \end{proof}
Note that as promised, the constructed condensation is ''nearly ordered''. It is now easy to describe the $CH$-construction.

\begin{Theorem}\Marginnote{Putting it all together}
Assume \CH\ and let $p_1,  p_2$ be non-isomorphic, selective ultrafilters containing $\min[FU(\mathbf{s})]$ and $\max[FU(\mathbf{s})]$ respectively.

Then there exists a stable union ultrafilter $u$ with $FU(\mathbf{s}) \in u$, $\min(u)=p_1$, $\max(u)=p_2$ and such that every ordered $\mathbf{t} \sqsubseteq \mathbf{s}$ has $FU(\mathbf{t}) \notin u$.
\end{Theorem}

\begin{proof} \plist
\item Assuming \CH, fix $\hl{(X_\alpha)_{\alpha < \omega_1}}$, an enumeration of $\mathp(\matf)$.

We argue by transfinite induction on $\beta < \omega_1$.

\item Assume for $\beta < \omega_1$ there are $\hl{ (FU(\mathbf{t}^\alpha))_{\alpha < \beta}}$  such that for all $\gamma < \alpha < \beta$
\begin{align*}
& \mathbf{t}^\alpha \sqsubseteq \mathbf{s} \\
& \mathbf{t}^\alpha \sqsubseteq^* \mathbf{t}^\gamma  \\
& FU(\mathbf{t}^\alpha) \subseteq X_\alpha \ \vee \ FU(\mathbf{t}^\alpha) \cap X_\alpha = \nil \\
& \min[FU(\mathbf{t}^\alpha)] \in p_1 \wedge \max[FU(\mathbf{t}^\alpha)] \in p_2  \\
& FU(\mathbf{t}^\alpha)\ \mathbf{s}\mbox{-meshed}
\end{align*}

\item Pick a cofinal sequence $(\alpha(n))_{n\in \omega}$ in $\beta$.

\item Applying Lemma \ref{lem:unordered main lemma} to $X:= X_\beta$ and $(FU(\mathbf{t}^{\alpha(n)}))_{n\in \omega}$ there exists $\mathbf{t}^\beta$ sufficient to continue the induction.

\item It should not be difficult to check that the resulting sets will generate a union ultrafilter as desired.
\pliste \end{proof}

Finally is useful to realize that the choice of the sequence $\mathbf{s}$ is not all that special.

\begin{Corollary}[The main theorem] \Marginnote{ \normalsize The main theorem} \label{thm:unordered union}
Assume \CH. For any two non-isomorphic, selective ultrafilters $p_1, p_2$ there exists a stable union ultrafilter $u$ which is not ordered, such that $\min(u)=p_1$ and $\max(u) = p_2$.
\end{Corollary}

\begin{spoiler}
The preceding theorem can be applied after using Theorem \ref{thm:blass homogeneity} to make sure that there is an apropriate sequence.
\end{spoiler}

\begin{proof} \plist
\item To invoke the preceding theorem it is sufficient to generate a \hl{suitably meshed sequence $\mathbf{s}$} with $\min[FU(\mathbf{s})] \in p_1, \max[FU(\mathbf{s})] \in p_2$.
\item For this consider the analytic set
\begin{align*}
\{ X \subseteq \omega \ | \ (\exists \mathbf{s} )\ & X = \max[FU(\mathbf{s})] \cup \min[FU(\mathbf{s})]  \\
& \mbox{and $FU(\mathbf{s})$ is $\mathbf{s}$-meshed}  \}.
\end{align*}

\item By Theorem \ref{thm:blass homogeneity} there exist \hl{$f$-alternating $Y_i \in p_i$} ($i\in 2$) such that set of $f$-alternating subsets of $Y_1 \cup Y_2$ is either included or disjoint from the analytic set.

\item There exist $f$-alternating $Y_i \in p_i$ ($i\in 2 $) such that $Y_0 \cup Y_1$ lies in the analytic set.
\begin{prose}
The argument is just as in the proof of the main lemma, i.e, it suffices to check that for any $Y_i \in p_i$ ($i\in 2$)
$ \min\!^{-1}[Y_1] \cap \max\!^{-1}[Y_2] $ must include $FU(\mathbf{s})$ for some suitably meshed sequence.
\end{prose}
\plist
\item Let $Y_i \in p_i$ be as in the previous step.

\item Recall that $\bbar{ \min\!^{-1}(p_1)} \cap \bbar{\max\!^{-1}(p_2)}$ is a closed subsemigroup; so we can find an idempotent ultrafilter therein.

\item Therefore there exists $\hl{FU(\mathbf{v})} \subseteq \min^{-1}[Y_1] \cap \max^{-1}[Y_2]$ by the Galvin-Glazer Theorem \ref{thm:galvin-glazer}.

\item Then there exists a condensation of $v$ to an $f$-alternating, meshed sequence $\mathbf{s}$, i.e., with $FU(\mathbf{s})$ being $\mathbf{s}$-meshed (with respect to $\mathbf{x}$).
\plist
\item For the inductive step $n \in \omega$ assume that for $k=\sum_{i<n} i$ there are $(s_i)_{i< k}$ with increasingly meshed graphs of sizes $1$ through $n-1$.

Pick $2n$-many elements from $FU(\mathbf{v})$ as follows:

\item First pick $(v_{i_j})_{j<n}$ past everything so far and then pick $(v_{i_j})_{n-1<j<2n}$ past additionally the ones just chosen and define
\[
s_{k+j} := v_{i_j} \cup v_{i_j+n}.
\]
\item Then $s_k,\ldots, s_{k+n}$ is an $n$-witness.

\item By construction, $\min[ \mathbf{s} ] \cup \max[ \mathbf{s} ] $ is $f$-alternating.
\pliste
\item Therefore, the $f$-alternating subsets are never disjoint from the anaytic set.
\item By Remark \ref{rem:alternating sets in selectives}, we find $Y_0, Y_1$ as desired.
\pliste
\item This completes the proof.
\pliste \end{proof}

Andreas Blass suggested an alternative proof for this last corollary sketched below.

\begin{proof} \plist
\item Given selectives $p_1,p_2$ there exists a permutation of $\omega$ simultaneously mapping $p_i$ to $p_i'$ ($i \in 2$) with $\min[FU(\mathbf{s})] \in p_1'$ and $\max[FU(\mathbf{s})] \in p_2'$; $p_1', p_2' $ are again selective. (Here, $\mathbf{s}$ is the previously fixed sequence).

\item The main theorem now gives a suitable $u'$ for $p_1', p_2'$.

\item But the natural extension of the permutation to $\matf $ yields an additive isomorphism on $FU(\mathbf{s}) \in u'$ mapping $u'$ to a union ultrafilter $u$ with $\min(u)=p_1$ and $\max(u)=p_2$.

\item Since additive isomorphisms preserve all the desired properties, this completes the proof.
\pliste \end{proof}

We can modify our construction to yield the following.

\begin{Theorem}
There exists an unordered union ultrafilter that is stable via $\min$. In particular, stability via $\min$ does not imply orderedness of a union ultrafilter.
\end{Theorem}
\begin{proof} 
\plist
\item We can modify the proof of the main lemma in the spirit of the (first) proof of the last corollary; compare the footnote in the proof of the main lemma. 

\item That is, in the inductive step of the Galvin-Glazer argument first choose an ordered sequence (of length $n+1$) followed by an $\mathbf{s}$-meshed witness (of length $n$) past this sequence. Finally add the elements from the ordered sequence to the $\mathbf{s}$-meshed witness just as in the proof of the corollary.

\item The ultrafilter $u$ resulting from this modified construction is of course still stable; in particular, it has the Ramsey property for pairs.
\item To show that it is stable via $\min$, let $f \in \mathbb{F}^\omega, g\in \omega^\omega$ with $f(s) < g\circ \min(s)$ on a set in $u$; for simplicity, we may assume that this set is all of $\mathbb{F}$.
\item Consider $\{ (s<t) \ | \ f(s\cup t) = f(s) \}$
\item Then there exists $A \in u$ with $A^2_<$ included in this set.
\plist
\item By the Ramsey property for pairs, we get a homogeneous set $A$. \item Fix $s\in A$
\item All $t\in \sigma(s)$ have $f(s\cup t ) < g\circ \min(s \cup t) = g\circ \min(s)$.
\item So on some $B\in u$, $f(s \cup t)$ is constant.
\item But now for any $t< t'$ in $B$, we get
\[
f( (s\cup t) \cup t') = f( s \cup (t \cup t') ) = f(s \cup t).
\]
\item In other words, $(s \cup t, t')$ is in the above set.
\pliste
\item So for ordered pairs, the value of $f$ on $A$ only depends on $\min$.
\item By construction of $u$, we find $\alpha < \omega_1$ such that $FU(\mathbf{s}^\alpha) \subseteq A$.
\item Then $f(s)$ depends only on $\min(s)$ on $FU(\mathbf{s}^{\alpha+1} )$.
\plist
\item Check that due to the modified construction every element of $\mathbf{s}^{\alpha +1}$ is a union of elements in $\mathbf{s}^{\alpha }$ where the first part is ordered with respect to the other parts.
\item Hence the value of $f$ depends only on that first part, i.e., only on $\min$.
\pliste
\pliste 
\end{proof}

As promised the assumption of the continuum hypothesis can be weakened.

\begin{Remark}
Dropping the prescribed selective ultrafilters in Lemma \ref{lem:unordered main lemma}, the modified consequent can be derived using Cohen forcing in the form of finite condensations of $\mathbf{s}$; using Lemma \ref{lem:preimages are meshed} it is not difficult to do some additional bookkeeping to ensure that the $\min$-image and the $\max$-image of the constructed union ultrafilter will be selective. 

Therefore, the above kinds of union ultrafilters already exist assuming $cov(\mathscr{M}) = \mathfrak{c} $ alone, in particular under weak versions of Martin's Axiom. For a detailed argument very much like the sketch we just proposed see \cite[Theorem 5]{BlassHindman87}.
\end{Remark}

To conclude this final section, we state some questions that remain open.

\begin{Question}

\plist
\item It is known that $\min$ and $\max$ of a union ultrafilter are not-near coherent $P$-points, the $\max$-image is rapid, cf.~\cite[Theorem 38]{Blass20092581}. Given such $P$-points on $\omega$, does there (say under \CH) exist a union ultrafilter mapping to them via $\min$ and $\max$?
\item More vaguely, do stronger assumptions hold for $\min$ and $\max$?
\item Most importantly, do there exist union ultrafilters that are not stable?
\item Is stability via $\min$ equivalent to stability?
\item Does the canonical partition property imply orderedness?
\pliste
\end{Question}
The first and second question are obviously related. A partition theorem for $P$-points similar to Theorem \ref{thm:blass homogeneity} can be found in \cite{Blass87-1} and strengthened, cf.~\cite[Theorem 4.10]{Krautzberger09}. This might be helpful in attacking the first question, especially if something can be improved regarding the second question. The third question seems to be an entirely different beast. It is much more difficult since the Galvin-Glazer Theorem so easily helps to construct almost condensation just as was done in the main result. It would seem to require a somewhat new proof of Hindman's Theorem to tackle stability.

\section*{Acknowledgments}

This work evolved out of the author's Ph.D.~thesis \cite{Krautzberger09} written under the supervision of Sabine Koppelberg at the Freie Universität Berlin and supported by the NaF\"oG grant of the state of Berlin. The author is also very grateful for the support of Andreas Blass especially during a visit to the University of Michigan, Ann Arbor, in the winter 2007/2008 with the support of the DAAD. The author also wishes to thank the organizing committee of BLAST 2010 for their financial support to attend the conference.

\bibliographystyle{alpha}


\end{document}